\documentclass[11pt,letterpaper]{article}
\usepackage[hmargin=1in, vmargin=1in]{geometry}
\usepackage{graphicx}
\usepackage{amsmath}
\usepackage{dsfont}
\usepackage{indentfirst}
\usepackage{sidecap}
\usepackage{float}
\usepackage{multicol}
\usepackage{caption}
\usepackage{subcaption}
\usepackage[hyphens]{url}
\usepackage{hyperref}

\usepackage{fancyhdr}
\title{3D Printing for Math Professors and Their Students}
\author{Edward Aboufadel \\ Grand Valley State University \and Sylvanna V. Krawczyk \\ California State University, Sacramento \and Melissa Sherman-Bennett \\Bard University at Simon's Rock}
\date{Version 1.1, October 2013}

\begin{document}
\maketitle

\section{Introduction}

Over the past few years, the 3D printing industry has experienced remarkable innovation, making it possible for anyone with access to a computer, design software, and a good idea to become a manufacturer. Recently, relatively inexpensive 3D printers have come on the market and are being purchased by colleges, universities, and libraries.  In this primer, we will describe a number of projects that can be completed with a 3D printer, particularly by mathematics professors and their students.  For many of the projects, we will utilize \emph{Mathematica}{\textsuperscript{\textregistered} to design objects that mathematicians may be interested in printing.

In general terms, 3D printing is an additive manufacturing technology. Traditional manufacturing methods start off with a block of material and whittle away unnecessary material to produce an object. 3D printing instead starts with nothing and builds up the object layer by layer, using only the material needed for the object itself. This saves on material costs and it allows for intricate designs that would be very difficult to accomplish with typical manufacturing methods.

The most affordable 3D printer models are extrusion-based models, where a thin plastic filament is fed into a heating element called an \emph{extruder} which melts the filament and lays down a thin trail of plastic onto the build plate. The extruder moves in the \emph{xy}-plane parallel to the build plate as specified by a file read by the printer. Once the extruder has deposited one layer of plastic onto the build plate, the build plate moves down slightly in the \emph{z} direction and the extruder lays a new horizontal layer on top of the previously extruded material. In this way, the printer stacks hundreds of horizontal layers on top of each other to manufacture user-designed objects.

In this primer, we will focus on projects designed for use with an extrusion-based 3D printer such as the Makerbot Replicator 2. The Replicator 2 has recently become affordable to mass markets, and because of this it is one of the most commonly available models when getting started with 3D printing.

\section{3D Printer Basics}

Before we describe the projects, it is important to explain some of the key aspects of using a 3D printer.

\vspace{12pt}

\noindent\textbf{Meshes:}  Computers define 3D objects through a mesh, which is a collection of polygons in $\mathds{R}^3$ that share edges, making one continuous surface. Both the type of polygon and number of them on the surface can be altered to increase or decrease the resolution of the object. To 3D print an object, the mesh must be saved as a STereoLithography (.stl) file, which is a file format that handles triangle meshes. The vertices and normal vectors of each triangle are recorded in an .stl file, and the file format is either binary (computer-readable only) or ASCII (human-readable).  These design files are shared on the Internet, in the same way that .pdf files are shared, at sites such as the Thingiverse (see Appendix A for a sample of sites).  An picture of a mesh can be seen at the beginning of Project 5 below.

\vspace{12pt}

\noindent\textbf{Slicing:}  After the object's mesh is produced, the .stl file must be processed by an intermediary program that converts our collection of triangles to a collection of horizontal layers, which the printer can then print. There are numerous ``slicers", among them the open-source \emph{Skeinforge} and MakerBot-developed \emph{MakerWare Slicer}, usually used in conjunction with an interface program such as \emph{ReplicatorG} with the former and \emph{MakerWare} with the latter. The interface allows you to preview and make small alterations to your model (scaling, rotating, etc.) before printing. The slicers produce G-code, which amounts to instructions to the 3D printer on where to move the extruder and how to extrude. \emph{ReplicatorG} allows the user to edit the G-code, technically giving the user the ability to manipulate how their model is printed. However, the nature of G-code makes it unfeasible to do significant alteration of printing in this manner, as it lists the extruder location, temperature, and speed for the entire printing process. From the G-code, the interface programs export an .x3g file, the native file format of the Replicator 2, which can be read by the printer.  (Other 3D printers utilize G-code differently.)  \emph{Skeinforge} and \emph{ReplicatorG} allow for more user interaction in the slicing and printing process, while \emph{Makerware} is more restrictive.

\vspace{12pt}

\noindent\textbf{\emph{Skeinforge} Profiles:}  \emph{Skeinforge} profiles are templates that govern general aspects of printing, such as maximum extruder temperature, extruder travel speed and layer thickness. They are used for all \emph{Skeinforge} slicings (all G-code generations in \emph{ReplicatorG} and high-quality slicing in \emph{MakerWare}). They are a more efficient way to customize how a model is printed, as the user edits and applies the profile of their choice to a model. The features are well-documented in the \emph{Skeinforge} wiki \cite {skeinforge}, but some key features are as follows (see figures below):

\vspace{12pt}

\begin{itemize}

\item \textbf{Raft} - produces a specified number of layers as a base for the model. This serves to attach models firmly to the build plate and minimizes damage to the model itself when removing it from the build plate. This is especially important for hollow models, or when just a surface is being printed.

\item \textbf{Supports} - produces external supports for the object wherever the slicer detects an overhang that would otherwise fall.

\item \textbf{Fill} - determines the solidity of a model. The G-code only specifies the perimeter of an object; this feature allows the user to control the percentage of the internal volume filled and in what manner.

\item \textbf{Skeinlayer} - produces a layer-by-layer diagram or animation of the extruder path as part of the slicing process, which is essentially a preview of how the model will print. This is especially important given the amount of time it takes to print most objects, as it allows the user to catch problems with the extruder path and avoid a failed print of the object. \url{gcode.ws} is an excellent online alternative to Skeinlayer if you are not using \emph{Skeinforge}, or if you want to view the extruder path without reslicing the object.
\end{itemize}

\begin{figure}[H]
        \centering
        \begin{subfigure}[b]{0.3\textwidth}
                \centering
                \includegraphics[width=\textwidth]{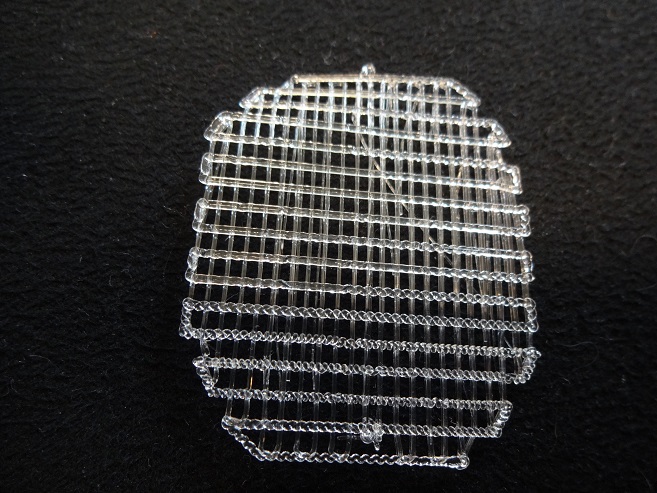}
        \end{subfigure}~
        \begin{subfigure}[b]{0.25\textwidth}
                \centering
                \includegraphics[width=\textwidth]{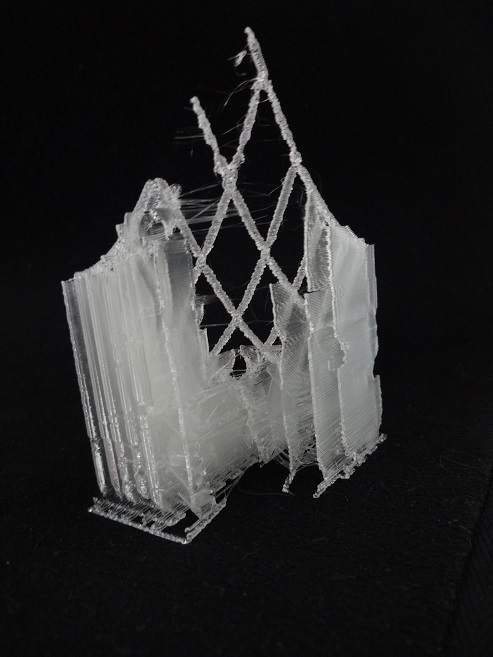}
        \end{subfigure}~
        \begin{subfigure}[b]{0.3\textwidth}
                \centering
                \includegraphics[width=\textwidth]{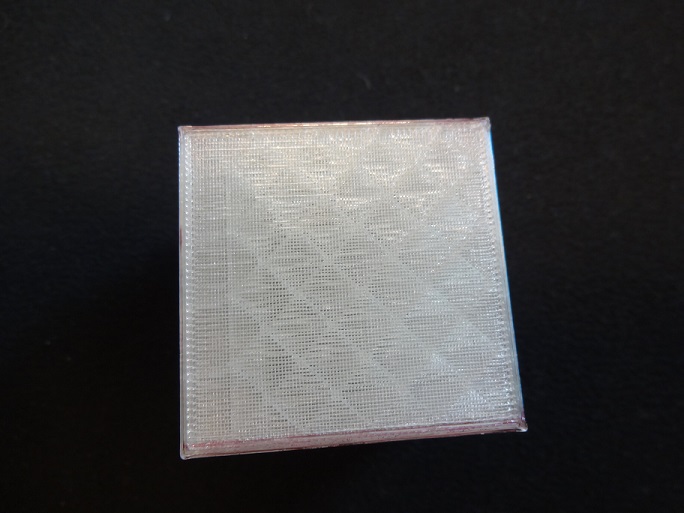}
        \end{subfigure}
\end{figure}

\centerline{Examples of a raft, support, and fill.}

\vspace{12pt}

\noindent\textbf{Printing:}  Once the .x3g file is produced, it can be printed directly from your computer if it's connected to the printer by a USB cable, or copied onto an SD card and then printed. The user's options once the print has started are fairly limited. We can pause the print, or cancel the print, but editing the object is not possible mid-print. The Replicator 2 will extrude plastic on to the build plate layer-by-layer until the object is completed or outside interference occurs.  For the projects described in this primer, typical printing times range from 15 minutes to one hour.  Because there is always a chance that an object will not print properly, it is important to keep an eye on the printing process while the 3D printer is working.

\section{Top 10 Tips for Successful 3D Prints}

Many of these tips will make more sense after trying out a few 3D prints.

\begin{enumerate}
\item Each horizontal slice (level curve) of your object must be a closed curve.  The overall print time of a model with non-closed level curves is far longer and the overall quality of the printed object is at risk if the object can be printed at all.
\item Temperature is key to proper filament flow.
\item Slicing takes a long time -- the simpler the mesh, the shorter your wait. Don't be alarmed if large or complicated objects take over a half hour to slice. However, if a relatively simple object, with simple level curves, is slicing very slowly, there's probably something wrong with your mesh.
\item Fine details tend to cause the extruder to jam.
\item It never hurts to take a look at the G-code. Make sure that rafts and supports are included in the G-code \emph{before} printing.
\item Smart design increases your chances of success (fewer supports, less drama with messy extrusion, less likely to fall over).
\item It helps to let filament ``drain" between prints. Any crust, goop, tangles, knots, etc. in the extruder can often be pushed out by loading the filament until a consistent stream of plastic flows through without catching on itself.
\item Make sure to load the spool in the proper orientation! The filament should unwind from the bottom of the spool into the guide tube; this maintains even tension on the filament when the extruder feeds in plastic. Loading in the other orientation introduces tangles on the spool and can inhibit the feed mechanism from pulling in enough plastic to continue the print.
\item Learn from your mistakes. If an object didn't print once, it's unlikely that a second try will go better. Alter the object (size, orientation, etc.), make some changes to the slicing profile and reslice, or perhaps reload the filament before printing it again.
\item If you're stuck, you're not alone -- check out the forums for help! (See Appendix A.)
\end{enumerate}

With these basic considerations in mind, the best way to familiarize yourself with the process of 3D printing is to jump right in with your first project. This primer to mathematical 3D printing is intentionally organized by project. Each section is designed to introduce a different method for generating mathematically derived objects from different types of data.

\section{Project 0:  Printing an .stl File Found Online}

If you have a new 3D printer, it probably came with previously-created .stl files that you can print.  Design files can also be found online at a variety of sites (see Appendix A).  For the Replicator 2, it is necessary to applying a slicing program to the .stl file (as described above) to create G-code, and then to convert the G-code to an .x3g file, followed by printing.  For this initial project, go through the process of printing a pre-made .stl file, following the directions with your printer.  Then we can move on to creating and printing our own objects.

\section{Project 1: Printing an Algebraically Defined Surface}

\includegraphics[width=\textwidth]{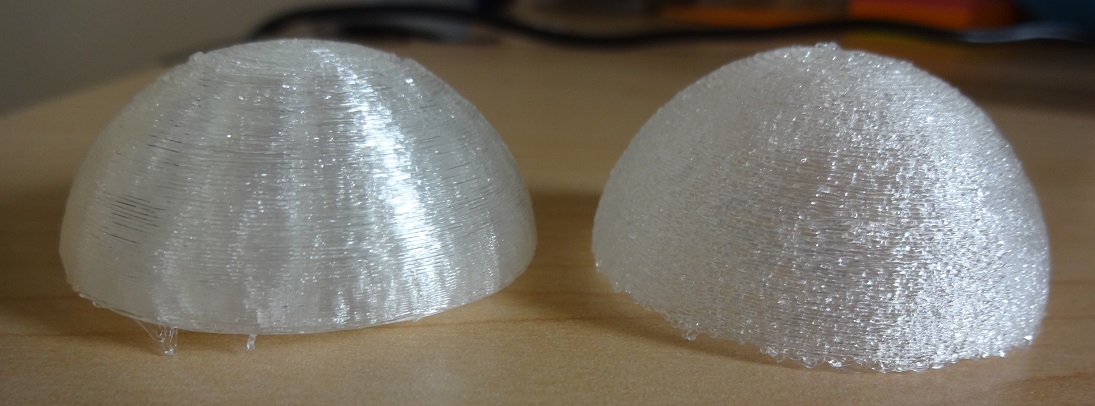}

We will start with the raw \emph{Mathematica} code to generate the mesh for the upper half of a sphere. To generate a 3D mesh for printing, it is sufficient to define a surface function, plot it with proper bounds, and export the plot as an .stl mesh:

\begin{verbatim}
1 surface =
   Plot3D[{Sqrt[1 - x^2 - y^2]}, {x, -1, 1}, {y, -1, 1},
   RegionFunction -> Function[{x, y, z}, x^2 + y^2 <= 1],
   BoxRatios -> Automatic]
2 Export["C:\\Users\\Owner\\Desktop\\hemisphere.stl", surface]}
\end{verbatim}

In our sample code above\footnote{Do not include the line numbers when using this or any other code in this paper!}:

\vspace{12pt}

\noindent \textbf{Line 1:} Designates a name for our surface (in this case, we have conveniently chosen to name it \texttt{surface}). \texttt{Plot3D} generates the surface of a sphere with radius 1.  \texttt{RegionFunction} limits the plot to the part of the sphere above the $x$-$y$ plane, bounding the bottom of the (now) hemisphere at $z=0$. This ensures a smooth base for the hemisphere.  \texttt{BoxRatios}, when set to Automatic, maintains the proportions of the hemisphere so it does not appear distorted.

\noindent \textbf{Line 2:} \texttt{Export} takes the plot of our hemisphere and generates an .stl file of the displayed plot. The filepath here places the file hemisphere.stl on the user's desktop. The ability to create .stl files is built into \emph{Mathematica}.

\vspace{12pt}

This algorithm produces the hollow shell of our hemisphere. To create a solid hemisphere, we can choose to add fill when slicing the .stl file we have produced.  For the hollow shell, if you have trouble printing the top of the hemisphere, try reducing the layer thickness when slicing, and increasing extruder temperature. You may also need to make the hemisphere smaller.

This is a very simple example of one type of algebraic surface. To print a different surface, we can alter the formula in line 2 to represent whatever three-dimensional region we wish. It is important to truncate the surface where it intersects the \emph{xy}-plane or some other horizontal plane, to provide a flat base for the printed surface. When designing the surface to print, it is beneficial to orient the surface so that the printer will create a stable, supported structure as it prints from the bottom up. For example, if we wanted to print a paraboloid with this method, it would be a wise decision to print the surface concave down in order to support the subsequent layers!

In terms of scale, this particular code is generating a hemisphere with radius 1. All slicing programs assume dimensions are in terms of millimeters, so this will correspond to a radius of 1 millimeter. Resize the hemisphere as you wish either in the slicing program or in the Mathematica code above.

Our method for producing 3D prints of algebraic surfaces was informed by Henry Segerman's method in his 3D printing article describing his approach for producing grid-like representations of algebraically-defined surfaces \cite{segerman}.

\section{Project 2: Image to Height}

In this project, we take an image and generate height information from the pixel values (i.e. how light or dark a pixel is, with 1 being white and 0 being black). These heights are plotted to create a mesh which is exported as an .stl file, as in our first project. The following works best with a very simple images with large areas of a single color, such as a cartoon character.  In the following code, we use the image ``picture.jpg".

\begin{verbatim}
1 image = Import["C:\\Users\\Owner\\Desktop\\picture.jpg"] // Binarize;
2 size = Import["C:\\Users\\Owner\\Desktop\\picture.jpg", "ImageSize"];
3 data = Table[
   10*ImageData[image][[i, j]], {i, 1, size[[2]], 1}, {j, 1,
    size[[1]], 1}];
4 object = ListPlot3D[data, BoxRatios -> Automatic]
5 Export["C:\\Users\\Owner\\Desktop\\picture.stl",object]
\end{verbatim}

You should see a 3D version of the simple picture you chose.
In our sample code above:

\vspace{12pt}

\noindent \textbf{Line 1:} The picture is imported from the desktop is designated as \texttt{image}. The \texttt{Binarize} option replaces all of the pixel values in the image above a certain threshold with the value 1, and all other values are replaced with 0. Using this technique, the image is rendered as closed regions of black and white.

\noindent \textbf{Line 2:} The dimensions of the picture are saved as a 2-element array \texttt{\{width, height\}}, which is called in \textbf{Line 3} as endpoints for a \texttt{Table} evaluation.

\noindent \textbf{Line 3:} \texttt{ImageData}[...][[i, j]] returns the pixel value of the pixel in the $i^{th}$ row and $j^{th}$ column of pixels in the picture, where the origin is the top left corner of the image. Because the image has been binarized, it outputs either 1 or 0. The \texttt{Table} function creates a subdivided array of the binarized values of every pixel. The parameter $i$ is evaluated first, and represents the $y$-coordinate of the pixel, not the $x$-coordinate. For each $i$ value, all of the $j$ values are evaluated (i.e. for each row, the pixel values of all of the columns are found), so each row of the array corresponds to a row in the picture.

\noindent \textbf{Lines 4-5:} The mesh is plotted, and the plot is exported in \textbf{Line 5}.  \texttt{BoxRatios} maintains the proportions of the original image so it does not appear distorted in the mesh.

\vspace{12pt}
\begin{figure}[H]

        \centering
        \begin{subfigure}[b]{0.3\textwidth}
                \centering
                \includegraphics[width=\textwidth]{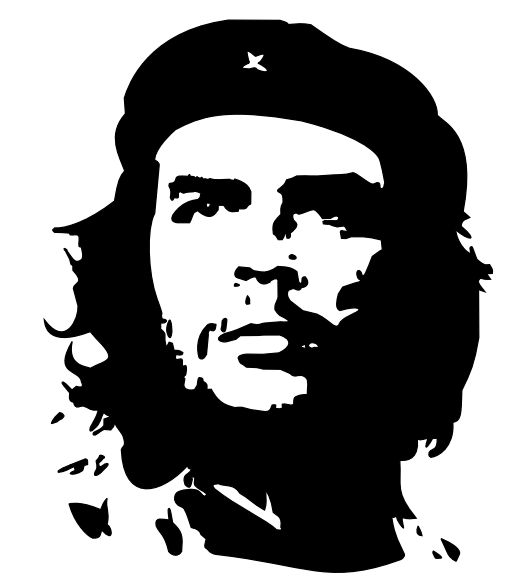}
        \end{subfigure}~
        \begin{subfigure}[b]{0.3\textwidth}
                \centering
                \includegraphics[width=\textwidth]{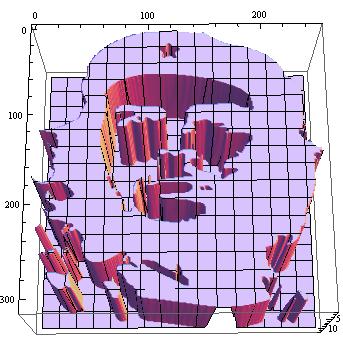}
        \end{subfigure}~
        \begin{subfigure}[b]{0.3\textwidth}
                \centering
                \includegraphics[width=\textwidth]{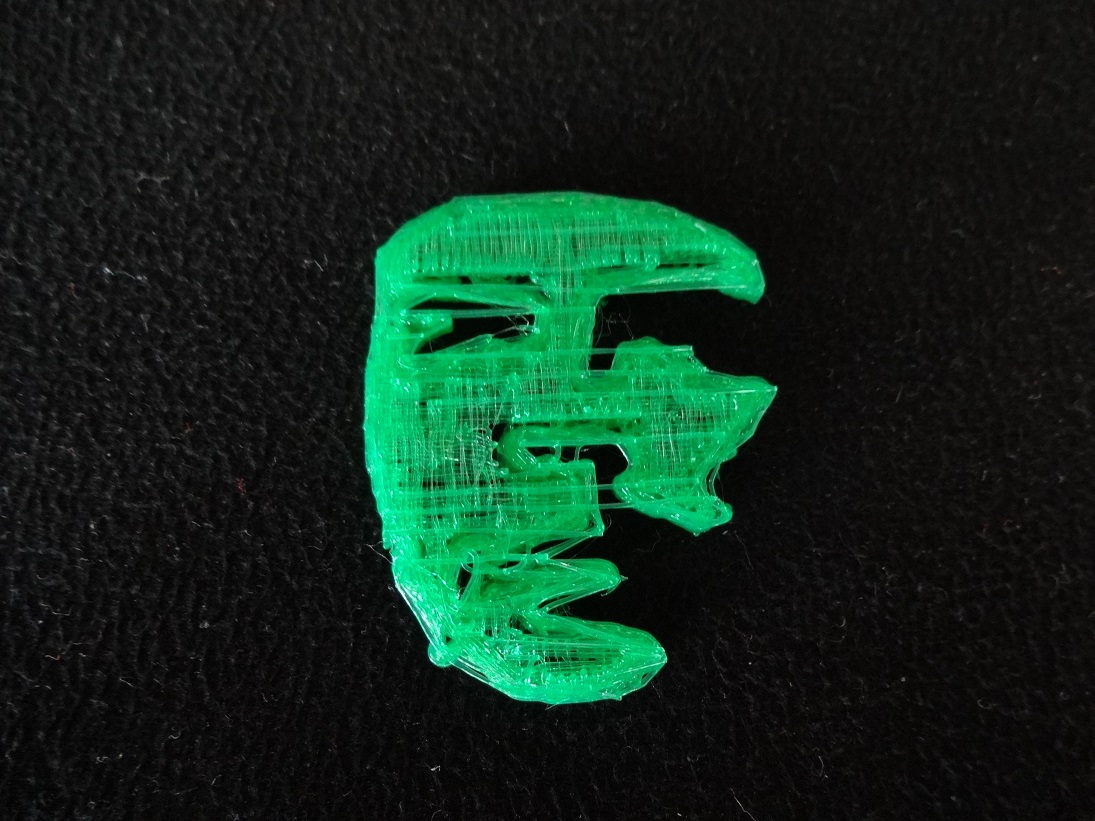}
        \end{subfigure}

\end{figure}

If the image you have chosen to work with is simple enough, the mesh you get from this process will have a few large plateaus representing your image. However, if you notice any small islands of elevated pixels in the mesh, you should be aware that these islands can be troublesome for the Replicator 2 to print. If you find that your mesh contains a large number of spikes (which would result from a complex or noisy image), the mesh may be too difficult to print. This is why choosing images with large swaths of continuous color is beneficial for printing binarized images. It is also possible to produce a mesh with multiple levels of height by using \texttt{ClusteringComponents} instead of \texttt{Binarize} (if \texttt{ClusteringComponents} is used, the image must be converted to grayscale before the function \texttt{ImageData} can be used). However, many meshes produced in this way are too detailed for successful printing.

\section{Project 3: Elevation Data}

A particularly interesting source of real life data which can be used to create a print is elevation data from terrain maps. Using a Google applet\footnote{\url{http://www.zonums.com/gmaps/terrain.php?action=sample}} developed by Zonum Solutions, we can extract the elevations at specified latitudes and longitudes, which we can use to print a geographic terrain with the Replicator 2.

In the elevation applet, we choose the range of latitudes and longitudes spanning the particular terrain we wish to capture. It works best to sample the data in a uniform grid, using at least 25 columns and rows. 50 rows and columns works quite well to achieve the resolution shown in the model of Mount St.~Helens below. Once the applet samples the appropriate points, we choose the option to view the data as a .txt file which contains all of the (latitude, longitude, elevation) coordinate triples for the chosen region.

\includegraphics[width=\textwidth]{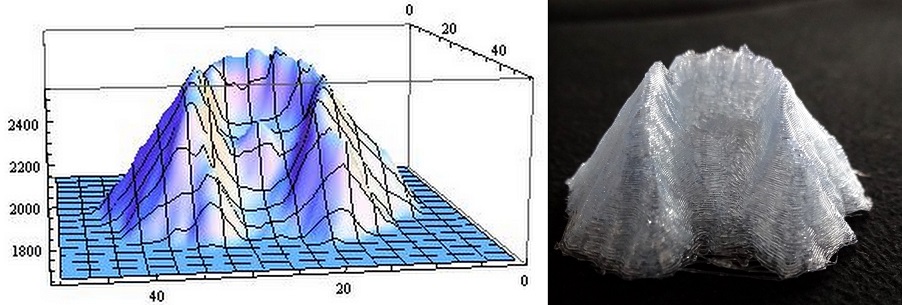}

Once we have the .txt file, we paste the data into \emph{Microsoft Excel} for convenient formatting. Using Excel's data import wizard, we want each element of the coordinate triples listed in its own cell starting with the first column and first row of the spreadsheet.

We should note that the data in our terrain file is organized in a very specific manner. The first latitude coordinate in our grid is held constant while the longitudes vary over the number of rows we have chosen. For example, if we choose to sample 30 rows by 45 columns, the first 45 rows of our .xls file will have the same latitude coordinate while the longitudes vary over the 45 sample points in our sampling grid. The next 45 rows of our .xls file will contain the latitude corresponding to the next row from our 30 row sample set. Thus, our .xls file will have $45\times30=1350$ rows total with 3 elements in each row corresponding to all of the latitude, longitude, and elevation points.  For our purposes, we only need the column containing the elevations, so we delete the first two columns of data in the spreadsheet. Once that has been accomplished, we save it as an .xls file using a memorable name -- we will be calling it in \emph{Mathematica}. For the purposes of this project, we will save it as the generic ``terrain.xls."

The arrangement of our elevation data will important for our 3D-plotting algorithm in \emph{Mathematica}:

\begin{verbatim}
1 input1 = Import["C:\\Users\\Owner\\Desktop\\terrain.xls"];
2 lengthX=50  (*insert number of coordinate triples with the same x-coordinate.
          In other words, the number of columns chosen in the applet.*)
3 size=2500  (*insert number of coordinate triples. In other words, the number of
          rows times the number of columns from your sample set.*)
4 minelevation=1500  (*insert the smallest value from your data set. In other words,
          the lowest elevation produced in the sample set.*)
5 data1 = Flatten[input1, 2];
6 data2 = ArrayPad[
   Table[(1/50)*
     Table[data1[[i]], {i, 1 + j, lengthX + j, 1}], {j, 0, (size)-(lengthX),
    lengthX}], 2, mineleveation];
7 terrain = ListPlot3D[data2];
8 plot = Show[terrain]
9 Export["C:\\Users\\Owner\\Desktop\\terrain.stl", plot]
\end{verbatim}

\noindent \textbf{Line 1:} The .xls file is imported containing the elevation data from the applet; the data is designated as \texttt{input1}.

\noindent \textbf{Line 2:} The number of repetitions in our data set is established. For our elevation data, the $x$-coordinate repeats until all of the longitudes have been sampled. Thus, the number of repetitions will be equal to the number of columns we chose to include in our grid in the applet.

\noindent \textbf{Line 3:} The total number of sampling points in our data set is established. For this project, this is the number of elevations the applet produced from our sample grid, which is the number of rows times the number of columns.

\noindent \textbf{Line 4:} The lower bound for our elevation dataset is defined. This becomes important for controlling how our plot will behave in \textbf{Line 6}.

\noindent \textbf{Line 5:} The \texttt{Flatten} command erases the levels within the array imported by our .xls file. When \emph{Mathematica} reads in this data, the array is represented as a nested list. For example, an array represented by $\bigl(\begin{smallmatrix} a & b & c \\ x & y & z \end{smallmatrix} \bigr)$ would be read into \emph{Mathematica} as $\{\{a, b, c\},\{x, y, z\}\}$. Each set of nested braces in this data list is a depth. The data from our .xls file is imported as a list of depth 3, with each elevation value being the only element in a 1 element list and all of the elevation values grouped as the first element of a 1 element list, and looks something like this: $\{\{\{a\}, \{b\}, \{c\}\}\}$. Our goal is to remove the innermost sets of braces so that we are left with a sequence of numbers of depth 1. The \texttt{Flatten} command with a second argument of 2 removes 2 levels of depth, giving us the simple list that we want.

\noindent \textbf{Line 6:} Starting with the innermost command, \texttt{Table} generates a table of information drawing from the list of elevations that we created in \textbf{Line 5}, with $\frac{1}{50}$ acting as a scaling factor for the elevations, which can be in the thousands.
This table will contain all of our elevation data, arranged in $i$ rows and $j$ columns to correspond to the original latitude-longitude grid.. The first set of braces contains the information needed to establish the $i$ rows. For each row, it grabs \texttt{lengthX} elements from our list, which represents all of the elevations at a given latitude. To decide where to begin the next row, the table references the second set of braces with the $j$ index. $j$ ranges from 0 to \texttt{size-lengthX} in increments of \texttt{lengthX}. \texttt{size-lengthX} represents the last elevation in our list that begins a new latitude and it indicates the element at which the last row of our plot should begin. $j$ runs in increments of \texttt{lengthX} because elements at the same latitude occur in our list as groups of \texttt{lengthX}, so after lengthX data points, a new row, corresponding to a new latitude, begins.  Outside of the \texttt{Table} command, we have the \texttt{ArrayPad} command. This borders the data table we have just created with a margin of rows and columns, the elements of which are all equivalent to the smallest elevation in our dataset. This creates a flat border around the mesh, which ensures that horizontal slices of our three-dimensional plot will be closed curves and that the Replicator 2 will cleanly print the mesh we are generating.

\noindent \textbf{Lines 7-9:} \texttt{ListPlot3D} takes the table we arranged in \texttt{Line 6} and plots it in $\mathds{R}^3$. For convenience, we name this \texttt{terrain}.  The plot is then exported to an .stl file called ``terrain.stl".

\vspace{12pt}

You may notice that our terrain mesh contains all of the desired features but the scales are slightly distorted. This is because \emph{Mathematica} plots the data symmetrically over the $xy$-plane, even if we did not choose our latitudes to cover an equal amount of land as our longitudes. The simplest way to avoid this is to choose square regions of terrain in the applet with equal numbers of columns and rows in our sample grid. Alternatively, if we do not wish to restrict ourselves to square plots, we can choose the number of columns and rows in our sample grid proportional to the ratios between the range of latitude and longitude. So, if we want to print a non-distorted segment of the Rocky Mountains covering three times more latitude than longitude, we should choose our grid to have three times as many rows as columns to maintain the proper proportions.

\section{Project 4: Kinect Portraits}
\vspace{12pt}
\begin{figure}[h]
\includegraphics[width=\textwidth]{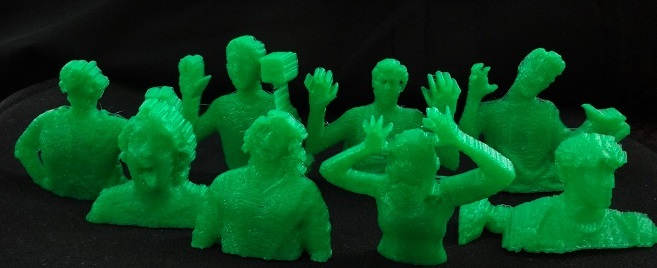}
\caption*{Examples of busts created using the process described here.}
\end{figure}
\vspace{12pt}

Another source for 3D data is the Microsoft Kinect\texttrademark camera system, which is an accessory for the Microsoft XBox\texttrademark.  The Kinect uses an infrared emitter and sensor to register how far away an object is \cite{kinectsensor1, kinectsensor2}. In addition to the standard color camera video feed, the Kinect's infrared depth feed can be accessed and depth data can be extracted from it.  The depth data can be plotted as a mesh in \emph{Mathematica} using the ideas in Project 3. The depth values from the Kinect are actual depths, so the meshes that we produce from the Kinect are essentially reliefs of the captured scene.

It takes a bit of work to enable a computer to extract data from a Kinect camera.  First, it is necessary to acquire a USB adapter/power source for the Kinect, which can be purchased through many online sources.  Second, a scripting program called \emph{Processing} needs to be installed on your computer, as well as drivers to allow the computer to communicate with the Kinect.  Plug-and-play does not work \cite{kinectpc}.

To be specific, the following steps work on a PC-based computer, although they take a bit of time to install:

\noindent\textbf{Installation of Drivers and \emph{Open NI}:}  An executable file that will install the needed Kinect drivers, as well as a software developer kit called \emph{Open NI}, can be found at the Brekel web site, at this URL: \url{http://www.brekel.com/kinect-3d-scanner/download}.  Select the link to the executable file in the ``All-in-one OpenNI Kinect Auto Driver Installer'' section.

\noindent\textbf{Installation of \emph{Processing} Scripting Program:}  This is an open-source program that we run to capture the Kinect depth data.  To install the program:  \url{http://processing.org/tutorials/gettingstarted/}.  We used the 32-bit version.  After getting Processing set up, when the file ``processing.exe'' is executed, a default folder ``Processing'' is created in My Documents.

\noindent\textbf{Installation of \emph{Processing} Helper Files:} A \emph{Processing} library called ``SimpleOpenNI'' needs to be included in the ``libraries'' subdirectory of the ``Processing'' directory.  The installation files can be found here:  \url{http://code.google.com/p/simple-openni/wiki/Installation}.  This library needs to be installed in ``Processing'' in My Documents, subdirectory ``libraries''.

Once all these files and drivers are installed, run \emph{Processing} and create the following code\footnote{partially based on code by from code.google.com/p/simple-openni/ and by Constantinos Miltiados of Architecture School/NTUA (Athens).}.  Save the code to a new directory as a .pde file:

\begin{verbatim}
1 import SimpleOpenNI.*;
2 SimpleOpenNI kinect;
3 PrintWriter output;
4 int kStep=3;
5 String myfile="Portrait.txt";

6 void setup()
7 {
8  kinect = new SimpleOpenNI(this);
9  kinect.enableDepth();
10 kinect.enableRGB();
11 background(200,0,0);
12 size(kinect.depthWidth() + kinect.rgbWidth() + 10, kinect.rgbHeight());
13 output = createWriter(myfile);
14 }

15 void draw()
16 {
17 kinect.update();
18 image(kinect.depthImage(),0,0);
19 image(kinect.rgbImage(),kinect.depthWidth() + 10,0);
20 }

21 void keyPressed()
22 {
23 output.println("Width "+kinect.depthWidth()/kStep);
24 output.println("Height "+kinect.depthHeight()/kStep);
25 output.println();
26 int index;
27 for (int i=0;i<floor(kinect.depthWidth()/kStep);i++) {
28 for (int j=0;j<floor(kinect.depthHeight()/kStep);j++) {
29      output.print(kinect.depthMap()[index]+"  ");
30      }
31     output.println();
32  }
33  output.println("[0,0,0]");
34  output.flush();
35  output.close();
36  println("Width "+kinect.depthWidth());
37  println("Height "+kinect.depthHeight());
38  println(kStep);
39  println("End of Program");
40  exit();
41 }
\end{verbatim}

\noindent \textbf{Lines 1-5:}  These are some initial commands to tell \emph{Processing} which helper files to use in order to receive data from the Kinect.  In \textbf{Line 4}, we specify we will extract the depth of every \texttt{kStep}\textsuperscript{th} pixel of the video stream, in order to save processing speed.  ``Portrait.txt'' is the data file that will be saved to the same directory as the .pde code file.

\noindent \textbf{Lines 6-14:}  When these commands are run, communication is opened between the computer and the Kinect.  Both cameras on the Kinect are enabled.

\noindent \textbf{Lines 15-20:}  These commands will create two camera images on the computer screen -- one from the typical RGB camera, the other from the infrared depth camera.

\noindent \textbf{Lines 21-32:}  In order to capture a depth data ``snapshot'', the user presses any key while the program is running.  When a key is pressed, the code in these lines is activated, which sends the depth data to ``Portrait.txt'', and terminates the program.  An item of note is that there are web sites indicating that the raw data exported here needs to be converted to a real distance, but the formula given on those sites is incorrect.  It appears that one person posted a faulty formula to a discussion board, and this formula was copied to several other sites. No conversion is necessary, because the depth data is automatically exported in terms of millimeters from the Kinect cameras.  Lines 23 and 24 are used to write the width and height of the ``depth picture'' to the output file.

\noindent \textbf{Lines 33-41:}  These commands are necessary to save the data file, followed by the printing of some debugging information to the computer screen.  At the beginning of the text file, the width and height of the sampled data will be printed, which tells you how many rows and columns of data you're dealing with. This information is important to know when importing the data to Mathematica.

\vspace{12pt}

Once the data is saved to a .txt file, it should be copied and pasted into \emph{Microsoft Excel} and saved as an .xls or .xlsx file.  Then, the following \emph{Mathematica} code can be used to convert the data to an .stl file, in a manner similar to Project 3:

\begin{verbatim}
1 input1 = Import["C:\\Users\\Owner\\Desktop\\Portrait.xls"];
2 bound[x_] := If[1500 >= x >= 20, x, 1500];
3 data1 = Flatten[input1, 2];
4 data2 = ArrayPad[
   Table[Table[(1/10) bound[ data1[[i]]], {i, 1 + j, 160 + j,
      1}], {j, 0, 33920, 160}], 2, 150];
5 bust = ListPlot3D[data2]
7 Export["C:\\Users\\Owner\\Desktop\\Portrait.stl", bust]
\end{verbatim}

\noindent\textbf{Line 2:} This line creates a function that will be used to limit the span of our data. If the input value isn't in the range specified here, it will be replaced by the value after the 2\textsuperscript{nd} comma; if it is within the specified range, it will be left alone. We here use a range of 20 to 1500, but this will need to be altered for your specific data. Play with these until the range encompasses the entire subject you wish to print and you have the desired level of detail.

\noindent\textbf{Line 4:} This line is the same as Line 6 in Project 3 (see said project for detailed explanation), with the exception of the application of the \texttt{bound} function. This function is simply applied to the depth data and acts as described above. The range for the $i$ parameter is the width of the depth data, while the range of the $j$ parameter is the (height*width) $-$ width, just as with Project 4. Note that the value for the margin created by \texttt{ArrayPad} is $\frac{\text{upper bound}}{10}$, the same as the highest depth value. This closes any open curves at the edge of the data plot.

\vspace{12pt}

Things to note about this procedure:
\begin{itemize}
\item The final print is the mirror image of the original, because of how the data plots in Mathematica. If desired, use a reflection transformation to ``flip'' the image before creating the .stl file.
\item The Kinect has trouble rendering shiny things, such as eyeglasses.
\item Don't cross anything in front of your body because it will create occlusions that cut into the print.
\item To make limiting the data easy, have your subject in front of a flat wall or sheet.
\item Larger prints have more potential for picking up detail; remember that the scale in \emph{Mathematica} is in millimeters.
\item Don't be discouraged if you see little detail at first. The data must sometimes be limited severely before detail can be seen, depending on the range of depth values you've captured. Try limiting only one end of the data (i.e. having only an upper bound), or reducing the range of data to 500 and go from there.
\end{itemize}

\section{Project 5: Write Your Own STL file}

\begin{center}
\includegraphics[width=0.5\textwidth]{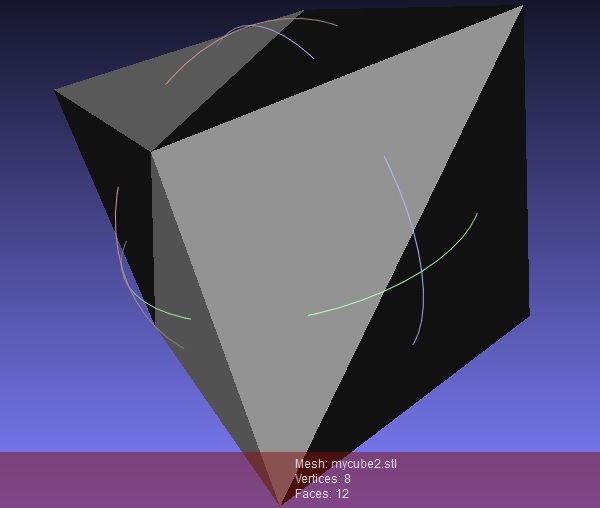}
\end{center}

An .stl file contains data that describes triangles in $\mathds{R}^3$, by indicating the three vertices of each triangle, as well as a unit outer normal vector.  While \emph{Mathematica} and other programs can generate .stl files automatically, usually in binary format, these files can also be created in a text editor in ASCII format.  For example, the following .stl code describes 12 triangles that correspond to six faces of the cube pictured above:

\begin{verbatim}





        solid mycube

          facet normal 0.0 -1.0 0.0
            outer loop
              vertex 0.0 0.0 0.0
              vertex 40.0 0.0 0.0
              vertex 0.0 0.0 40.0
            endloop
          endfacet

          facet normal 0.0 -1.0 0.0
            outer loop
              vertex 40.0 0.0 40.0
              vertex 40.0 0.0 0.0
              vertex 0.0 0.0 40.0
            endloop
          endfacet

          facet normal 0.0 1.0 0.0
            outer loop
              vertex 0.0 40.0 0.0
              vertex 40.0 40.0 0.0
              vertex 0.0 40.0 40.0
            endloop
          endfacet

          facet normal 0.0 1.0 0.0
            outer loop
              vertex 40.0 40.0 40.0
              vertex 40.0 40.0 0.0
              vertex 0.0 40.0 40.0
            endloop
          endfacet

          facet normal -1.0 0.0 0.0
            outer loop
              vertex 0.0 0.0 0.0
              vertex 0.0 40.0 0.0
              vertex 0.0 0.0 40.0
            endloop
          endfacet

          facet normal -1.0 0.0 0.0
            outer loop
              vertex 0.0 40.0 40.0
              vertex 0.0 40.0 0.0
              vertex 0.0 0.0 40.0
            endloop
          endfacet

          facet normal 1.0 0.0 0.0
            outer loop
              vertex 40.0 0.0 0.0
              vertex 40.0 40.0 0.0
              vertex 40.0 0.0 40.0
            endloop
          endfacet

          facet normal 1.0 0.0 0.0
            outer loop
              vertex 40.0 40.0 40.0
              vertex 40.0 40.0 0.0
              vertex 40.0 0.0 40.0
            endloop
          endfacet

          facet normal 0.0 0.0 -1.0
            outer loop
              vertex 0.0 0.0 0.0
              vertex 40.0 0.0 0.0
              vertex 0.0 40.0 0.0
            endloop
          endfacet

          facet normal 0.0 0.0 -1.0
            outer loop
              vertex 40.0 40.0 0.0
              vertex 40.0 0.0 0.0
              vertex 0.0 40.0 0.0
            endloop
          endfacet

          facet normal 0.0 0.0 1.0
            outer loop
              vertex 0.0 0.0 40.0
              vertex 40.0 0.0 40.0
              vertex 0.0 40.0 40.0
            endloop
          endfacet

          facet normal 0.0 0.0 1.0
            outer loop
              vertex 40.0 40.0 40.0
              vertex 40.0 0.0 40.0
              vertex 0.0 40.0 40.0
            endloop
          endfacet

        endsolid
\end{verbatim}

Each \texttt{facet} block of code describes one of the triangles, starting with the outer normal vector, and then the $(x, y, z)$ coordinates of the three vertices.  The above code can be saved as a file such as ``cube.stl'', then sliced and printed, creating a cube.

By writing code to manipulate string variables and export them to a text file, we can write our own .stl files using \emph{Mathematica}.  The following code provides us more flexibility than the native .stl exporter in \emph{Mathematica}.

\begin{verbatim}
1 myfile =  "C:\\Users\\Owner\\Desktop\\myobject.stl"

2 f [x_, y_] := Max[0, (15 - x^2 - y^2)*Exp[-(x/5)^2 - (y/5)^2] + 5]
3 xmin = -10.0;
4 xmax = 10.0;
5 ymin = -10.0;
6 ymax = 10.0;
7 res = 5.0;    (*higher number yields a finer mesh *)

8 startsolid = "solid Default";
9 endsolid = "endsolid Default";
10 startloop = "\n    outer loop";
11 endloop = "\n    endloop";
12 startfacet = "\n  facet normal ";
13 endfacet = "\n  endfacet";
14  avertex = "\n      vertex ";

15 nv[x1_, y1_, x2_, y2_, x3_, y3_] :=
      Normalize[Cross[{x2 - x1, y2 - y1, f[x2, y2] - f[x1, y1]},
        {x3 - x1, y3 - y1, f[x3, y3] - f[x1, y1]}]]

16 facet[x1_, y1_, x2_, y2_, x3_, y3_] :=
 StringJoin[startfacet,
  ToString[FortranForm[nv[x1, y1, x2, y2, x3, y3][[1]]]], " ",
  ToString[FortranForm[nv[x1, y1, x2, y2, x3, y3][[2]]]], " ",
  ToString[FortranForm[nv[x1, y1, x2, y2, x3, y3][[3]]]],  startloop,
  avertex, ToString[FortranForm[x1]], " ", ToString[FortranForm[y1]],
  " ", ToString[FortranForm[f[x1, y1]]], avertex,
  ToString[FortranForm[x2]], " ", ToString[FortranForm[y2]], " ",
  ToString[FortranForm[f[x2, y2]]], avertex,
  ToString[FortranForm[x3]], " ", ToString[FortranForm[y3]], " ",
  ToString[FortranForm[f[x3, y3]]], endloop, endfacet]

17 f1 = StringJoin[Table[facet[i, j, i + 1/res, j, i, j + 1/res],
      {i, xmin, xmax,1/res}, {j, ymin, ymax, 1/res}]];

18 f2 = StringJoin[Table[facet[i + 1/res, j + 1/res, i + 1/res, j, i, j + 1/res],
      {i, xmin, xmax, 1/res}, {j, ymin, ymax, 1/res}]];

19 allfacets = StringJoin[f1, f2];

20 stllines = Array["", 3];
21 stllines[[1]] = startsolid;
22 stllines[[2]] = allfacets;
23 stllines[[3]] = endsolid;

24 Export[myfile, stllines, "TEXT"]
\end{verbatim}

\noindent \textbf{Line 1:} The name of the .stl file is defined.

\noindent \textbf{Lines 2-7:} This is where to enter an algebraically-defined function whose surface we wish to print, as well as the function's domain.  In order for the slicer to recognize the facets, it is important to include decimal points to force floating point arithmetic.  The resolution of the mesh is set in \textbf{Line 7}.

\noindent \textbf{Lines 8-14:} In this part of the code, string variables are defined that will be used in the .stl file.  It is critically important that the exact number of blank spaces in the code are used (such as the four spaces in \textbf{Line 10}), or else the .stl file will not be interpreted properly.

\noindent \textbf{Line 15:}  This is a function to compute the unit normal vector, using the cross product from multivariable calculus.  It does not seem to matter whether this vector is pointing ``out'' or ``in''.

\noindent \textbf{Line 16:}  This code is used to define a function which will generate one of the facets, based on the $x$ and $y$ coordinates of the three vertices of the triangle.

\noindent \textbf{Lines 17-19:}  Here we have a looping through a lattice of $x$ and $y$ coordinates, combining string variables together that represent facets.  Text for all of the facets are stored in \texttt{allfacets}.

\noindent \textbf{Lines 20-24:}  In this part of the code, the preamble and end of the .stl file are combined with the code for the facets, and all of the .stl design code is then written to the text file in \textbf{Line 1}.

\vspace{12pt}

As with the earlier projects, it is important to pay close attention to where the object intersects the $xy$-plane. This code can be modified to produce .stl files for all of the previous projects, though only the algebraic surface version is given here. Note that though this process gives you more flexibility, it produces a larger .stl file than would be produced by Mathematica's export function, which can be problematic.

\section{Project 6: Printing Objects From Photos}

During Summer 2013, the authors developed a method to extract an object's depth information from a set of photographs taken from various perspectives around the object.  With this method, it is possible to 3D print an object based on photographs.  This method will be described in a subsequent paper which the authors are writing.  See \url{sites.google.com/site/aboufadelreu/} for updates.

\section{Other Mathematical Connections to 3D Printing}

Motivated mathematics faculty and students may be interested in pursuing the following questions:

\begin{enumerate}
\item A \emph{Jordan curve} is a simple closed curve, and according the Jordan curve theorem, every Jordan curve divides the plane into two pieces -- an interior and an exterior.  When an object represented in an .stl file is sliced, will each slice be a Jordan curve?

\item The G-code that is used for printing describes the path that the extruder takes in laying down the plastic used for the print.  To what extent is that path optimized in the sense of distance traveled, or through some other measure?  This question may be connected to the \emph{Traveling Salesman Problem}.

\item 3D printers are being used to create complex, mathematically-based puzzles akin to the Rubik's cube.  Can you design your own unique puzzle?  To learn more: \url{http://gigaom.com/2013/07/29/the-puzzle-masters-how-3d-printing-is-enabling-the-most-complex-puzzles-ever-created/}

\item If you Google \texttt{mathematics of 3D printing}, one of the links that is found is a blog written by ``Andy''.  He describes a conversation \cite{andy} with Prof. Alan Branford, Director of Studies in Mathematics and Statistics at Flinders University in South Australia, regarding 3D printing.  Andy asked Prof. Branford if there is a proof that ``we are able to express any 3D Printed object as layers of 2D planes''.  The response was: ``Yes, Fubini's Theorem''.  Is this correct?
\end{enumerate}











\section{Appendix A: Helpful URLs for 3D printing}
\begin{itemize}
\item \emph{Skeinforge} wiki:  \url{http://wiki.bitsfrombytes.com/index.php/Skeinforge}
\item Online G-code viewer:  \url{http://gcode.ws/}
\item Online .stl viewer:    \url{http://www.marcoantonini.eu/doku.php?id=3dmodel:stlviewer}
\item Thingiverse .stl repository:  \url{http://www.thingiverse.com/}
\item Shapeways 3D parts .stl repository:  \url{http://www.shapeways.com/3d_parts_database}
\item 3D Printing Forums:  \url{http://www.3dprinting-forums.com/}, \url{http://www.shapeways.com/forum/}, \url{http://support.makerbot.com/forums}, \url{http://groups.google.com/forum/#!forum/makerbot}
\end{itemize}

\section{Appendix B:  A Brief Selection of Interesting Mathematically-Connected 3D Printing Links}

Francesco de Comité's homepage:  \url{http://www.lifl.fr/~decomite/tdprinting/tdprinting.html}

Henry Segerman's page:  \url{http://www.ms.unimelb.edu.au/~segerman/}

Oliver Knill (of Harvard U.)'s pages on 3D printing:  \url{http://www.math.harvard.edu/~knill/3dprinter/exhibit.html}

Daniel Walsh's Geometric Curiosities:  \url{http://danielwalsh.tumblr.com/post/8350627231/printing-the-impossible-in-3d}

George Hart's 3D Fractals:  \url{https://www.simonsfoundation.org/multimedia/3-d-printing-of-mathematical-models/}

3D Printing with Blender:  \url{http://www.blendernation.com/2013/03/27/blender-2-67-feature-3d-printing-toolbox/}

Bathsheba Grossman:  \url{http://boingboing.net/2013/05/20/profile-of-math-inspired-3d-pr.html}

\section{Appendix C: Further Information}

For questions or comments about this paper, contact Prof. Edward Aboufadel of Grand Valley State University: \url{mailto:aboufade@gvsu.edu}.  This paper is available at \url{sites.google.com/site/aboufadelreu/}.  Also available are some of the codes described in this paper, without the line numbers.

This work was partially supported by National Science Foundation under Grant Number DMS-1262342, which funds a Research Experience for Undergraduates program at Grand Valley State University.  Any opinions, findings, and conclusions or recommendations expressed in this material are those of the author(s) and do not necessarily reflect the views of the National Science Foundation.

\end{document}